\documentclass{amsart}

\usepackage{amsthm,amssymb,latexsym,amsmath}
\usepackage{mathrsfs}
\usepackage{graphicx}
\usepackage{color}
\usepackage[brazil]{babel}
\usepackage[latin1]{inputenc}
\usepackage{graphicx}%
\usepackage{multirow}%
\usepackage{amsmath,amssymb,amsfonts}%
\usepackage{amsthm}%
\usepackage{physics}
\usepackage{mathrsfs}%
\usepackage[title]{appendix}%
\usepackage{xcolor}%
\usepackage{textcomp}%
\usepackage{manyfoot}%
\usepackage{booktabs}%
\usepackage{algorithm}%
\usepackage{algorithmicx}%
\usepackage{algpseudocode}%
\usepackage{listings}%
\usepackage{amsthm,amssymb,latexsym,amsmath}
\usepackage{mathrsfs}
\usepackage{graphicx}
\usepackage{color}
\usepackage[all]{xy}
\usepackage{graphicx}%
\usepackage{multirow}%
\usepackage{amsmath,amssymb,amsfonts}%
\usepackage{amsthm}%
\usepackage{mathrsfs}%
\usepackage[title]{appendix}%
\usepackage{xcolor}%
\usepackage{textcomp}%
\usepackage{manyfoot}%
\usepackage{booktabs}%
\usepackage{algorithm}%
\usepackage{algorithmicx}%
\usepackage{algpseudocode}%
\usepackage{listings}%

\newcommand{\CC}{{\mathbb C}}
\newcommand{\OO}{\mathcal O}

\newcommand{\fol}{\mathcal{F}}

\newcommand{\ZZ}{{\mathbb{Z}}}
\newcommand{\PP}{\mathbb{P}}

\newtheorem{lema}{Lemma}[section]
\newtheorem{cor}[lema]{Corollary}

\newtheorem*{teo1234}{Theorem}
\newtheorem*{teo1}{Theorem 1}
\newtheorem*{teo2}{Theorem 2}
\newtheorem*{teo3}{Theorem 3}
\newtheorem*{teo4}{Theorem 4}
\newtheorem*{teo5}{Theorem 5}

\theoremstyle{definition}
\newtheorem{remark}[lema]{Remark}
\newtheorem{defi}[lema]{Definition}
\newtheorem{exe}[lema]{Example}

\begin{document}

\title{
THE SCHWARTZ INDEX AND THE RESIDUE OF LOGARITHMIC FOLIATIONS ALONG A HYPERSURFACE WITH ISOLATED SINGULARITIES}

\hyphenation{lo-ga-ri-thmic}

\begin{abstract}
Given a compact complex manifold $X$, we prove a Baum-Bott type formula for one-dimensional holomorphic foliations on $X$ that are logarithmic along a hypersurface with isolated singularities. We show that the residues of these foliations can be expressed in terms of the Schwartz index of the vector fields that locally define them. Furthermore, in this context, we prove that the Schwartz index is positive when $\dim(X)$ is even and that the GSV index is positive when $\dim(X)$ is odd. As application, we show that the obstruction determined by the multiplicity of the isolated singularities of the invariant hypersurface, for the solution of Poincaré's problem in holomorphic foliations on $\PP^2$, is a more general fact, valid for holomorphic foliations defined on projective spaces of arbitrary even dimension. Additionally, we prove that the obstruction determined by the Euler characteristic for the existence of vector fields is even more comprehensive in the case of hypersurfaces with isolated singularities.
\end{abstract}

\author{Diogo da Silva Machado}

\address{\noindent  Diogo da Silva Machado\\
Departamento de Matem\'atica \\
Universidade Federal de Vi\c cosa\\
Avenida Peter Henry Rolfs, s/n - Campus Universitário \\
36570-900 Vi\c cosa- MG, Brazil} \email{diogo.machado@ufv.br}



\maketitle
\section{Introduction}

\hyphenation{ge-ne-ra-li-zes}
\hyphenation{cha-rac-te-ris-tic}

The index of a vector field at a singular point is one of the most fundamental concepts in the study of smooth manifolds, and its properties have important consequences in results such as the Poincaré-Hopf theorem. This theorem establishes a formula relating the total sum of the indices of a vector field to the Euler-Poincaré characteristic of the ambient manifold where the vector field is defined. It has given rise to a vast literature and it was the genesis of several lines of research, including the theory of vector fields on singular varieties and the Baum-Bott theory of residues, both of which are fundamental to the present work.

The study of vector fields on singular varieties was initiated by M.-H. Schwartz, who defined the \textit{Radial index} (also called the \textit{Schwartz index} in the literature) for a special class of vector fields and proved a version of the Poincaré-Hopf theorem for singular varieties (\cite{SSS1}, \cite{SSS2}). This concept was later generalized by various authors in different contexts (see, for example, \cite{AGG}, \cite{ABC00}, \cite{ABC}, \cite{Ebe}, \cite{Ebe2}, \cite{King}, \cite{SeaSuw2}). Other indices for vector fields on singular varieties also exist, such as the GSV index, introduced by X. Gómez-Mont, J. Seade, and A. Verjovsky \cite{GSV} for hypersurfaces with isolated singularities and extended in \cite{SeaSuw1} to complete intersections. The GSV index generalizes the Poincaré-Hopf index and preserves the important property of stability under perturbations of singularities. Many generalizations of the GSV index have been obtained by different authors (see, for example, \cite{a03}, \cite{a09}, \cite{cm3}, \cite{a08}, \cite{a07}).

On the other hand, Baum-Bott theory of residues, initiated by the pioneering work of P. Baum and R. Bott \cite{BB1}, localizes the Euler-Poincaré characteristic of the ambient manifold in the singular set of a holomorphic foliation defined on it. More generally, for a one-dimensional holomorphic foliation $\fol$ with isolated singularities on an
$n$-dimensional complex manifold $X$, the residue is obtained as localization of the normal sheaf of $\fol$, and the residue at each singular point $p \in Sing(\fol)$ is given by the Milnor number of $\fol$ at $p$ (i.e., the Milnor number at $p$ of a holomorphic vector field $v$ that defines $\fol$ locally at $p$). Furthermore, if $X$ is compact, one has

\begin{eqnarray}\label{bbformula}
  \displaystyle \hspace*{2.0 cm}\int_{X}c_n(TX - T\fol) \,\,\, =\,\,\, \sum_{p\in Sing(\fol)} \mu_{p}(\fol),
\end{eqnarray}

\noindent where the total Chern class is $c_{\ast}(TX - T\fol) = c_{\ast}(TX) \cdot c_{\ast}(T\fol)^{-1}$ and $\mu_{p}(\fol)$ denotes the Milnor number of $\fol$ at $p$. Recall that $\mu_{p}(\fol)$ also coincides with Poincaré-Hopf index at $p$ of a holomorphic vector field $v$ that defines $\fol$ locally at $p$ and can be expressed in terms of Grothendieck residue $$\mu_{p}(\fol) =  dim_{\CC}\frac{\OO_{n,p}}{\left\langle a_1,\ldots, a_n\right\rangle}=\mbox{Res}_p\left[\begin{array}{cccc} \det(Jv)\\ a_1,\ldots,a_n\end{array}\right ],$$
\noindent where $(a_1,\ldots,a_n)$ denotes local coordinate of $v$ near $p$ and $Jv$ is its Jacobian matrix.

\hyphenation{res-tric-ted}

Subsequently, Baum-Bott theory of residues was developed for foliations restricted to their invariant variety (see for example \cite{ds}, \cite{Suw2}), where each residue --- also referred to as a {\it relative residue} --- is obtained as the localization of appropriate vector bundles or virtual bundles.  These residues often coincides with one of the indices of the vector fields that locally define the foliation, such as the GSV index. In the global viewpoint, the sum of these residues (and these indices) equals the integration over invariant variety of the top Chern class of the (virtual) bundles in question, in the spirit of the Baum-Bott formula (1).

In this work, we prove a Baum-Bott type formula for holomorphic foliations of dimension one on complex compact manifolds that are logarithmic along a hypersurface with isolated singularities. In this case, we have shown that the residue at a  singularity of foliation is expressed in terms of the Schwartz index of vector field that locally define the foliation. More precisely, we prove Theorem 1, which will be demonstrated in Section \ref{a1}:

\begin{teo1} Let $X$ be an $n$-dimensional complex manifold, $D\subset X$ a hypersurface with isolated singularities and $\fol$ a one-dimensional holomorphic foliation on $X$, with isolated singularities, which leaves $D$ invariant. If $\varphi$ is a homogeneous symmetric polynomial of degree $n$, then 
\begin{itemize}
\item [(I)]  For each $p\in Sing(\fol) \cup Sing(D)$ there exists the residue $\mbox{Res}_{\varphi}(\fol,D, p)$,   which is a complex number that depends only on $\varphi$ and on the local behavior of the leaves of $\fol$ near $p$.

\medskip

\item [(II)] If $X$ is compact, then
$$ 
\int_{X} \varphi (TX - [D] - T_{\fol}) = \sum_{p} \mbox{Res}_{\varphi}(\fol,D,  p) $$
\noindent where the sum is taken over all points $p$ of $Sing(\fol) \cup Sing(D)$ and $[D]$ denotes the line bundle of the divisor $D$.

\medskip

\item [(III)] In particular, if $X$ is compact and $\varphi = c_n$ then
{\small {
\begin{eqnarray}\nonumber
\hspace{1.7 cm} \int_{X} \hspace{-0.10 cm} c_n (TX - [D] - T_{\fol}) \hspace*{0.05cm} = \hspace*{-0.25cm}\sum_{p\in Sing(\fol)}\hspace{-0.35 cm} \mu_{p}(\fol)  - \hspace*{-0.25cm} \sum_{p\in  S(\fol, D)}\hspace{-0.35 cm} Sch_p(\fol,D)\, + \, (-1)^{n}\hspace{-0.35 cm}\sum_{p\in  Sing(D) } \hspace{-0.35 cm} \mu_p(D)
\end{eqnarray}
}}
\noindent where $Sch_p(\fol,D)$ denotes the Schwartz index of $\fol$ along $D$ in $p$, $\mu_p(D)$ is the Milnor number of $D$ at $p$ and $S(\fol, D) = (Sing(\fol) \cap D) \cup Sing(D)$.
\medskip
\item [(IV)] Also, one has
\begin{itemize}
\item [(a)] If $p\in Sing(\fol)\cap (X - D)$, then $$\mbox{Res}_{c_n}(\fol,D,  p) =  \mu_p(\fol);$$
\item [(b)] If $p\in Sing(\fol) \cap D_{reg}$, then
$$\mbox{Res}_{c_n}(\fol,D,  p) =  \mu_p(\fol) - Sch_p(\fol, D);$$ 
\item [(c)] If $p\in Sing(\fol) \cap Sing(D)$, then
$$\mbox{Res}_{c_n}(\fol,D,  p) =  \mu_p(\fol) - Sch_p(\fol, D) + (-1)^n\mu_p(D);$$ 
\item [(d)] If $p\in Sing(D)  - Sing(\fol)$, then
$$\mbox{Res}_{c_n}(\fol, D,  p) = - Sch_p(\fol, D) + (-1)^n\mu_p(D);$$.
\end{itemize}
\item [(V)] In particular, if $X=\PP^n$, with $deg(\fol) = d$ and $deg(D) = k$ then
\begin{eqnarray}\nonumber
\hspace{0.5cm}\displaystyle\sum_{p\in  S(\fol, D) } \hspace{-0.4cm} Sch_p(\fol, D) \,\,\,= \,\,\, \sum_{i=0}^{n-1}\left(1 - (1-k)^{n-i}\right)d^i \,\,\,\, + \,\,\,\,(-1)^n\hspace{-0.45 cm}\displaystyle\sum_{p\in Sing(D)}\mu_p(D). 
\end{eqnarray}
\end{itemize}
\end{teo1}

\begin{remark} \label{rema01}
Item (V) of Theorem 1 establishes a formula for the total sum of the Schwartz indices in terms of the degree of the foliation and the degree of the invariante hypersurface. By \cite[Ch.IV, Proposition 1.9]{Suw2} we have 
\begin{eqnarray}\label{form1}
Sch_p(\fol,D) = GSV_p(\fol,D) + (-1)^n\mu_p(D).
\end{eqnarray}
\noindent Thus we also obtain the following formula for the total sum of the of GSV indices
\begin{eqnarray}\label{lu}
\displaystyle\sum_{p\in  S(\fol, D) } GSV_p(\fol,D) \,\,\,= \,\,\, \sum_{i=0}^{n-1}\left(1 - (1-k)^{n-i}\right)d^i, 
\end{eqnarray}
\noindent which is a version of formula in \cite[Proposition 6.2]{DSM}.
\end{remark}

Considering singular holomorphic foliations, logarithmic along a hypersurface with isolated singularities, we proved that the Schwartz index is positive, if $n$ is even, and that the GSV is positive, if $n$ is odd. We refer to Subsection \ref{secschw} for details on the Schwartz index in the context of holomorphic foliations, and Subsection \ref{sec111}  for definition and details on the GSV index. More precisely, we prove the Theorem 2, which will be demonstrated in Section \ref{a4}:

\begin{teo2} \label{tt2} Let $X$ be an $n$-dimensional complex manifold, $D\subset X$ a hypersurface with isolated singularities and $\fol$ a one-dimensional holomorphic foliation on $X$, with isolated singularities, which leaves $D$ invariant. Then for all point $p\in Sing(\fol)\cap D$ one has
\begin{itemize}
\item [(i)]   $Sch_p(\fol, D)>0$, if $n$ is even; 
\bigskip
\item [(ii)] $GSV_p(\fol, D)>0$, if $n$ is odd.
\end{itemize}
\end{teo2}

\hyphenation{in-va-riant }
\hyphenation{mul-ti-pli-ci-ty}
\hyphenation{differential}
\hyphenation{lea-ving}

As applications of Theorem 1 and Theorem 2, we provide some results related to Poincaré's problem for foliations. We recall that the Poincar\'e's problem is the question of bounding the degree of a hypersurface $D \subset \PP^n$, invariant by a holomorphic foliation $\fol$ on $\PP^n$, in terms of the degree of the foliation. This problem originated with Poincaré \cite{PPP}, in his search for a characterization of algebraic differential equations that are algebraically integrable.

In the last three decades, we have witnessed the development of a vast literature on Poincaré's Problem.  For example, D. Cerveau and A. Lins Neto \cite{Alc} obtained the inequality
\begin{eqnarray}\label{eh}
deg(D) \leq deg(\fol) + 2,
\end{eqnarray}
for a one-dimensional foliation $\fol$, on the complex projective plane $\PP^2$, leaving invariant an algebraic curve $D$, with at most nodal singularities.

We observe that every nodal singularity is an isolated singularity with multiplicity 1. We prove that this condition determines an obstruction to the solution of the Poincaré problem and generalizes formula (\ref{eh}) for arbitrary even dimension. More precisely, as consequence of Theorem 1 and Theorem 2, we prove the following theorem: 

\hyphenation{sin-gu-la-ri-ties}
\begin{teo3} \label{antC}
Let $\fol$ be a foliation of dimension one on $\PP^n$, with isolated singularities, and $D \subset \PP^n$ a hypersurface with isolated singularities such that $\fol$ leaves $D$ invariant, with $Sing(D) \subset Sing(\fol)$. Suppose that $n$ is even and that all singular points of D have multiplicity 1. Then 
$$
deg(D) \leq deg(\fol) + 2.
$$
\end{teo3}

\hyphenation{associa-ted}

M. Soares \cite{Soa} considered, in addition to the degree, other characters associated to an invariant hypersurface and proved the following Poincaré inequality, for a holomorphic foliation in the complex projective plane:

\begin{teo1234}[M. Soares, \cite{Soa}]
Let $D \subset \PP^2$ be an irreducible hypersurface with isolated singularities and degree $\deg(D) = k$. Let $\fol$ be a foliation of dimension one on $\PP^2$, with isolated singularities and degree $\deg(\fol) = d$. If $\fol$ leaves $D$ invariant, then
\begin{eqnarray}\label{9}
k^2 - 2k \hspace*{0.25cm} - \sum_{p\in Sing(D)}(\mu_p(D)- 1) \hspace*{0.25cm}\leq \hspace*{0.25cm}\displaystyle kd. 
\end{eqnarray}
\end{teo1234}

\noindent In his demonstration M. Soares also considered $Sing(D) \subset Sing(\fol)$.  As applications of Theorem 1 and Theorem 2, we have a generalization of the formula (\ref{9}) for any even dimension:

\begin{teo4} Let $\fol$ be a foliation of dimension one in $\PP^n$, with isolated singularities and degree $\deg(\fol) = d$. Let $D \subset \PP^n$ be a hypersurface of degree $\deg(D) = k$ with isolated singularities and such that $\fol$ leaves $D$ invariant. Suppose $n$ is even and $Sing(D) \subset Sing(\fol)$. One has 

\begin{eqnarray} \label{oneo}
\hspace*{-0.0cm} \displaystyle \sum^{n-1}_{j=0}\hspace*{-0.10cm}\binom{n}{j}\hspace*{-0.05cm}(-k)^{n-j} \hspace*{-0.12cm} -  \hspace*{-0.45cm}\sum_{p\in Sing(D)}\hspace*{-0.5cm} \left(\mu_p(D) -1\right) \hspace*{-0.01cm}\leq \hspace*{-0.01cm}\displaystyle \sum^{n-1}_{i=1}\sum^{i-1}_{j=0}\binom{i}{j}(-1)^{i-j+1}k^{i-j}d^{n-i}.
\end{eqnarray}
\end{teo4}

Further results related to Poincaré's problem can be found, for example, in \cite{Brun},  \cite{Carni}, \cite{Vileh}, \cite{CE}, \cite{per}, \cite{mgSoa}. In \cite{DSM}, a result similar to Theorem 3 was obtained, considering some additional assumptions. Additionally, a characterization of formula (\ref{oneo}) in terms of the Schwartz and GSV indices was also obtained.

Still considering the consequences of Theorem 1 and Theorem 2, we obtained some applications to other problems related to holomorphic foliations. We know that the Poincaré-Hopf theorem states that the Euler-Poincaré characteristic of a compact oriented manifold is a measure of the obstruction to the construction of a tangent vector field to the manifold, without singularity. The positivity of the Poincaré-Hopf index in holomorphic vector fields gives the these obstruction a much broader aspect, which can be expressed as follows:

``{\it A compact complex manifold $X$ admitting a global holomorphic vector field with $s$ isolated zeros satisfies $\chi(X)>s$.''}

\noindent This more comprehensive aspect of the obstruction is the manifestation of the fact that there are very few global holomorphic vector fields on a compact complex manifold, according to a result due to J. Carrel and  D. Lierbermann \cite{rell}.

If $X$ is a singular complex variety we don't know if the techniques developed in \cite{rell} apply so directly, since in the singular context the complex of sheaves of regular holomorphic forms is not isomorphic to the (classical) Koszul complex. However, by a generalization of the Poincaré-Hopf theorem for vector fields on singular varieties the Euler-Poincaré characteristic represents the obstruction for the construction of global vector fields, without singularity (see for example \cite[Corollary 2.4.1]{BarSaeSuw}). 

As a consequence of Theorem 2, we obtain the following theorem that making more comprehensive these obstruction determined by the Euler-Poincaré characteristic:  

\begin{teo5}
Let $D$ be a complex compact hypersurface with isolated singularities $p_1,\ldots,p_{s_1}$ in a complex manifold $X$ of even dimension. If $D$ admits a holomorphic vector field singular at the $p_j's$ and possibly at some other smooth points $q_1,\ldots,q_{s_2}$ then
$$
\chi(D) > s_1+s_2.
$$
\end{teo5}

In Corollary \ref{coro}, we proved that the {\it Tjurina number} of an $\fol$-invariant hipersurface is a lower bound for the GSV index of $\fol$ and, in the projective case, we obtained a lower bound for GSV index in terms of the (local) multiplicity of $D$.

\section{Preliminaries}
\hyphenation{theo-ry}
\hyphenation{di-men-sio-nal}

\subsection{Singular one-dimensional holomorphic foliations }
\begin{defi} 
Let $X$ be a connected  complex manifold. A one-dimensional holomorphic foliation $\fol$ is
given by the following data:
\begin{itemize}
  \item[$i)$] an open covering $\mathcal{U}=\{U_{\alpha}\}$ of $X$;
  \item [$ii)$] for each $U_{\alpha}$ a holomorphic vector field $v_\alpha \in TX|_{U_{\alpha}}$;
  \item [$iii)$]for every non-empty intersection, $U_{\alpha}\cap U_{\beta} \neq \emptyset $, a
        holomorphic function $$f_{\alpha\beta} \in \mathcal{O}_X^*(U_\alpha\cap U_\beta);$$
\end{itemize}
such that $v_\alpha = f_{\alpha\beta}v_\beta$ in $U_\alpha\cap U_\beta$ and $f_{\alpha\beta}f_{\beta\gamma} = f_{\alpha\gamma}$ in $U_\alpha\cap U_\beta\cap U_\gamma$.
\end{defi}
We denote by $K_{\fol}$ the line bundle defined by the cocycle $\{f_{\alpha\beta}\}\in \mathrm{H}^1(X, \mathcal{O}^*)$. Thus, a one-dimensional holomorphic  foliation $\fol$ on $X$ induces  a global holomorphic section $s_{\fol}\in \mathrm{H}^0(X,T_X\otimes K_{\fol})$. The line bundle $T_{\fol}:= (K_{\fol})^* \hookrightarrow T_X$ is called the  tangente  bundle of $\fol$. The singular set of $\fol$ 
is $Sing(\fol)=\{s_{\fol}=0\}$. We will assume that $codim (Sing(\fol))\geq 2$.

Given a hypersurface $D \subset X$, we say that the foliation $\fol$ \emph{leaves} $D$ \emph{invariant} if it satisfies the following condition: for all $x\in D - Sing(D)$, the vector $v_{\alpha}(x)$ belongs to $T_xD$, with $x\in U_{\alpha}$. In \cite{MD},\cite{MD3}, \cite{DSM} the foliations that leaves $D$ invariant were also called {\it logarithmic foliations} along $D$.

A foliation on a complex projective space $\mathbb{P}^n$  is  called a projective foliation. Let $\fol $ be  a   {\it projective foliation}   with tangent bundle $T_{\fol}=\mathcal{O}_{\mathbb{P}^n}(r)$. The integer $d:=r+1$ is called the degree of $\fol$.

It is well known that if $\fol$ is a one-dimensional holomorphic foliation on $\PP^n$, with isolated singularities and degree $\deg(\fol) = d \geq 0$, the classical Baum-Bott formula says that the total sum of Milnor numbers of $\fol$ is given by (see for example \cite[Remark 5.1]{Soar})
\begin{eqnarray}\label{2121}
\sum_{p\in  Sing(\fol)} \mu_p(\fol) = d^n+d^{n-1}+\ldots + d+ 1.
\end{eqnarray}

\subsection{The Schwartz index} \label{secschw}
Let us recall the definition of the Schwartz index \cite[Ch.2, 2.1]{BarSaeSuw}. For simplicity, we consider the case where the singular space is a hypersurface $D$, with isolated singularities, on an n-dimensional complex manifold $X$. 

Given an isolated singularity $p\in Sing(D)$, let $(D,p)$ be the germ of $D$ at $p$ and assume that it is embedded in some complex space $\CC^n$. Let $v$ be the restriction to (a representative of) $D$ of a holomorphic vector field $\tilde{v}$ on a neighborhood of $p$ in $\CC^n$, which is tangent to regular part of $V$. If $v$ is everywhere transverse to the link $K = D \cap S_{\epsilon}$ of $p$ in $V$, its Schwartz index is 1. Otherwise, let $v_0$ be the restriction to (a representative of) $D$ of other holomorphic vector field $\tilde{v}_0$ on a neighborhood of $p$ in $\CC^n$, which is tangent to regular part of $V$ and suppose that when restricted to the link $K'=D \cap S_{\delta}$, $\delta < \epsilon$, $v_0$ is the unit outwards-pointing normal vector field on $K'$ in $V$. Let $H \subset D - \{p\}$ be the cylinder in D bounded by $K$ and $K'$ and let $\vartheta$ be the vector field on $H$ that  has finitely many singularities in the interior of $H$ and restricts to $v_0$ on $K'$ and $v$ on $K$.  The {\it Schwartz index} of $v$ at $p\in D$ is defined by
$$
Sch_p(v,D) = 1 + d(v_0,v)      
$$

\noindent where $d(v_0,v) \in \ZZ$ is the total Poincaré-Hopf index of $\vartheta$ in $H$, called {\it difference between} $v_0$ and $v$.

Note that if $p$ is a regular point of $D$, then $Sch_p(v,D)$ coincides with the Poincare-Hopf index.

\begin{remark}
The definition of the Schwartz index for the case of a variety with isolated singularities and continuous (or $C^{\infty}$) vector fields is made in an analogous way. Generalizations in different contexts  can be found in \cite{AGG}, \cite{ABC00}, \cite{ABC}, \cite{Ebe}, \cite{Ebe2}, \cite{King}, \cite{SeaSuw2}.
\end{remark}

Let $\fol$ be a foliation of dimension one on $X$, with isolated singularities leaving invariant the hypersurface $D$. If $\tilde{v}$ defines $\fol$ locally at $p$, according to \cite[Ch.IV, Corollary 7.15]{Suw2} the Schwartz index of $v=\tilde{v}|_{D}$ at $p\in D$ does not depend on the choice of such a vector field $\tilde{v}$. The {\it Schwartz index of $\fol$ in $p$ along $D$}, denoted by $Sch_{p}(\fol, D)$, is defined as $Sch_{p}(v, D)$.

\subsection{The GSV-index} \label{sec111}
X. Gomez-Mont, J. Seade and A. Verjovsky \cite{GSV} introduced the GSV-index for a holomorphic vector field over an analytic hypersurface with isolated singularities, on a complex manifold, generalizing the (classical) Poincar\'e-Hopf  index. The concept of GSV-index was extended to continuous vector fields on more general contexts. For example, J. Seade and T. Suwa in \cite{SeaSuw1}, defined the GSV-index for vector fields on isolated complete intersection singularity germs. J.-P. Brasselet, J. Seade and T. Suwa in \cite{a03}, extended the notion for vector fields defined in certain types of analytic subvariety with non-isolated singularities. 

\hyphenation{using}
In \cite{a08}   X. Gomez-Mont  defined the homological index of holomorphic vector fields on an analytic hypersurface with isolated singularities, which coincides with the GSV-index. There is also the virtual index, introduced by D. Lehmann, M. Soares and T. Suwa \cite{a07}, defined via Chern-Weil theory, that can be interpreted as the GSV-index. 
M. Brunella \cite{a09} also presents the GSV-index for foliations on complex surfaces  by a different approach.

\hyphenation{sin-gu-la-ri-ties}
\hyphenation{sin-gu-la-ri-ty}
\hyphenation{de-fi-ning}

Let us recall the  definition of the GSV-index  (\cite{BarSaeSuw}, Ch.3, 3.2). 
Let $D$ be a hypersurface  with isolated singularities  on an $n$-dimensional complex manifold $X$ and let $v$  be a holomorphic vector field on $X$ with isolated singularities, and logarithmic along $D$. Given a singular point $p_0 \in Sing(D)$, let $h $  be an analytic function  defining  $ D$ on a neighborhood $U_0$ of $p_0$. The gradient vector field $\overline{grad}\,(h)$ is nowhere vanishing away from $p_0$, because $p_0$ is an isolated singularity.  

Denote by $v_{\ast}$ the restriction of $v$ to the  regular part $D_{reg} =  D - Sing(D)$ of $D$. On the neighborhood $U_0$,  suppose that  the vector field $v$ is non-singular away from $p_0$. Since $v$ is logarithmic along $ D$, we have that $\overline{grad}\,(h)(z)$ and $v_{\ast}(z)$ are linearly independent at each point  $z\in U_0 \cap ( D - \{p_0\})$. Assume that $(z_1,\ldots,z_n)$ is a system of complex coordinates on $U_0$ and consider 
 $$S_{\varepsilon} = \{z= (z_1,\ldots,z_n):\,\,\mid\mid z-p_0\mid\mid \,\, = \,\, \varepsilon\}$$  
 the sphere   sufficiently small so that $K =  D \cap S_{\varepsilon} $ is the link of the singularity of $ D$ at $p_0$ (see, for example, \cite{JM}). It is a $(2n-1)$-dimensional real oriented manifold. By using the Gram-Schmidt process, if necessary, the vector field $v_{\ast}$ and $\overline{grad}\,(h)$  define a continuous map

$$
\phi_{v}:=(v_{\ast},\overline{grad}\,(h)): K \longrightarrow W_{2,n+1}
$$  

\noindent where $W_{2,n+1}$ is the Stiefel manifold of complex $2$-frames in $\CC^{n+1}$. The GSV-index of $v$ in $p_0 \in D$, denoted by $GSV_{p_0}(v, D)$, is defined as the degree of map $\phi_{v}$.

\begin{remark}
In the above definition $D$ has dimension greater than 1 or it has only one branch. When $dim(D) =1$ and $D$ has several branches the GSV index is defined as the Poincaré-Hopf index of an extesion of $v$ to a Milnor fiber (see \cite{Bru1}, \cite{Bru}, \cite{KS}). 
\end{remark}

Now, given a one-dimensional foliation $\fol$ on $X$, with isolated singularities, 
leaving invariant the hypersurface $D$, let $v$ be a holomorphic vector field defining $\fol$ on a neighborhood $U$ of $p_0$ in $X$. The $GSV_{p_0}(v, D)$ does not depend on the choice of such a vector field $v$, local representative of the foliation $\fol$ at $p$ (see \cite[Ch.IV, Corollary 7.15]{Suw2}). The GSV-index of $\fol$ in $p_0$ along $D$, denoted by $GSV_{p_0}(\fol, D)$, is defined as $GSV_{p_0}(v, D)$. Moreover, in this case, by \cite[Ch.IV, Theorem 7.16]{Suw2} we have
\begin{eqnarray}\label{foj}
\displaystyle\int_{D} c_{n-1}(TX - [D]- T_{\fol}) = \sum_{p\in  S(\fol, D) } GSV_p(\fol,D),
\end{eqnarray}
\noindent where $S(\fol, D) = (Sing(\fol) \cap D) \cup Sing(D)$. On the other hand, according to the (generalized) Adjunction Formula \cite[Theorem 5.6.3]{BarSaeSuw}, the Euler-Poincaré characteristic of $D$ is given by 
\begin{eqnarray}\label{gaf}
\chi(D) = \int_{D} c_{n-1}(TX - [D]) + (-1)^{n}\sum_{p\in Sing(D)}\mu_p(D),
\end{eqnarray}
\noindent and using (\ref{form1}) one has the following formula \cite[Ch.IV, Theorem 7.16]{Suw2}
{\small{
\begin{eqnarray}\label{fof}
\displaystyle\int_{D}\hspace{-0.1 cm} c_{n-1}(TX - [D]- T_{\fol}) -\int_{D}\hspace{-0.25 cm} c_{n-1}(TX - [D])  + \chi(D) \,\,\,\,\,\,= \hspace{-0.1 cm} \sum_{p\in  S(\fol, D) }\hspace{-0.4 cm} Sch_p(\fol,D)
\end{eqnarray}
}}

\section{Proof of Theorem 1}\label{a1}

Given  $p\in Sing(\fol) \cup Sing(D)$, let $U$ be an open neighborhood of $p$ in $X$ disjoint from the other point  of $Sing(\fol) \cup Sing(D)$. Letting $U_0 = U-\{p\}$ and $U_1 = U$, we consider the covering $\mathcal{U} = \{U_0,U_1\}$ of $U$ and the associated $\check{\mbox{C}}$ech-de Rham complex $A^{\bullet}(\mathcal{U})$. 

For each $k=0,1$, let  $\nabla_k, \nabla'_k$ and $\nabla''_k$ be connections for bundles $\displaystyle TX, [D]$ and $T_{\fol}$, respectively, on $U_k$. If we set $\nabla_k^{\bullet} = (\nabla''_k, \nabla'_k, \nabla_k)$ the characteristic class $\varphi (TX - [D] - T_{\fol})$  is represented by the cocycle 
$$\varphi(\nabla^{\bullet}_{\ast}) = (\varphi(\nabla^{\bullet}_{0}), \varphi(\nabla^{\bullet}_{1}), \varphi(\nabla^{\bullet}_{0}, \nabla^{\bullet}_{1} ))$$ 
\noindent in $\check{\mbox{C}}$ech-de Rham cohomology group
$H^{2n}(A^{\bullet}(\mathcal{U}))$.

On the other hand, by taking a smaller $U$, if necessary, we may assume that there is a generator $v$ of $\fol$ on $U$, where $v$ is a non-vanishing holomorphic vector field on $U_0$. 

\noindent We claim that $\eval{TX}_{U_0}, \eval{[D]}_{U_0}$ {\it and} $\eval{T_{\fol}}_{U_0}$ {\it are holomorphic $v$-bundles}. Of course,  the restrictions $\eval{TX}_{U_0}$ and $\eval{[D]}_{U_0}$ becomes holomorphic $v$-bundles, respectively, by the actions
$ \alpha_v: \Gamma(U_0,\eval{TX}_{U_0}) \longrightarrow \Gamma(U_0,\eval{TX}_{U_0})$,  defined by $\alpha_v(w) = [v,w]$, and $ \alpha'_v: \Gamma(U_0,\eval{[D]}_{U_0}) \longrightarrow \Gamma(U_0,\eval{[D]}_{U_0})$ defined as follows: let $f\in \OO_X(U)$ such that $f = 0$ is an equation for $D$, in U,  and consider the sequence
\begin{eqnarray} \label{seq3030}
\centerline{
\xymatrix{
0\ar[r]& \eval{Tf}_{U_0} \ar[r]&   \eval{TX}_{U_0}\ar[r]^{\text{\large{$\pi$}}}  &    \eval{[D]}_{U_0} \ar[r] &0,  \\
}}
\end{eqnarray}  
\noindent where $\eval{Tf}_{U_0}$ is a bundle of vectors tangent to the fibers of $f$. For $\vartheta \in \Gamma(U_0,\eval{[D]}_{U_0})$, we may write $\vartheta = \pi(\omega)$ with $\omega$ a section of $\eval{TX}_{U_0}$ and we set 
$\alpha'_v(\vartheta) = \pi([v,\omega])$. On the other hand, to prove that $\eval{T_{\fol}}_{U_0}$ is a $v$-bundle, we can assume that $\eval{T_{\fol}}_{U_0}$ is isomorphic to trivial line bundle $U_0\times \CC$ (by taking a smaller $U$, if necessary) and consider the action $\alpha''_v$ given by $\alpha''_v(h) = v(h)$ for $h\in \Gamma(U_0,\eval{T_{\fol}}_{U_0}) = \OO_X(U_0)$.  This establishes the claim.

Now, since $\eval{TX}_{U_0}, \eval{[D]}_{U_0}$ and $\eval{T_{\fol}}_{U_0}$ holomorphic $v$-bundles,  we can consider $\nabla_0, \nabla'_0$ and $\nabla''_0$ to be $v$-connections for $\displaystyle TX, [D]$ and $T_{\fol}$, respectively, on $U_0$. Therefore, it follows from Bott vanishing theorem (see for example \cite[Ch.II, Theorem 9.11]{Suw2}) that
$$
\varphi(\nabla^{\bullet}_{0}) \equiv 0.
$$

\noindent Consequently, the cocycle $\varphi(\nabla^{\bullet}_{\ast})$ is in relative  $\check{\mbox{C}}$ech-de Rham complex $A^{2n}(\mathcal{U}, U_{0})$ and it defines a class in the relative cohomology $H^{2n}(U, U - \{p\}; \CC)$, which we denote by $\varphi_{p}((TX - [D] - T_{\fol})|_{U},\fol)$. The residue $\mbox{Res}_{\varphi}(\fol,D,  p)$ is defined as the image of that class by Alexander duality $H^{2n}(U,U-\{p\};\CC)  \simeq H_0(\{p\};\CC) = \CC$ and the item (I) of theorem is proved.

On the other hand, if $X$ is compact, then the proof of item (II) follows from the \cite[Ch.II, Proposition 3.11]{Suw2}. Of course, in this case, just consider the commutative diagram

$$
\xymatrix{H^{2m}(X,X-S;\CC) \ar[r]^{j^{\ast}} \ar[d]^{A_{M_S}} &  H^{2m}(R,\partial R;\CC) \ar[d]^{L_{R}} \\H_{0}(S,\CC) \ar[r]^{i_{\ast}} & H_{0}(R,\CC)} 
$$
with $R = X$, $\partial R = \emptyset$, $S = Sing(\fol)\cup Sing(D)$ and the decompition $\displaystyle H_{0}(S,\CC) = \displaystyle \oplus_{p \in S} H_{0}(\{p\},\CC)$. In this case, we consider the inclusions $i:S\hookrightarrow R$ and $j:(R,\partial R) \hookrightarrow (M,M-S)$ and the dualities of Alexander ($A_{M_S}$) and Lefschetz ($L_{R}$).

Now, consider $X$ compact and $\varphi = c_n$. By using the definition of Chern class of virtual bundle, we obtain
\begin{eqnarray}\nonumber
\displaystyle\int_{X}c_{n}(TX - [D]- T_{\fol})
&=&\displaystyle\int_{X}\sum^{n}_{j=0}c_{n-j}(TX - [D])c_1( T_{\fol}^*)^j
= 
\\\nonumber &&\\ \nonumber &=&
\displaystyle\int_{X}\sum^{n}_{j=0}\left(\sum^{n -j}_{i=0}c_{n-j-i}(TX)c_1([D]^*)^i\right)c_1( T_{\fol}^*)^j.
\end{eqnarray}

\hyphenation{fol-lowing}
\hyphenation{trian-gu-lar}

\noindent Since $[D]$ and $T_{\fol}$ line bundles, we have $c_1([D]^*)c_1( T_{\fol}^*) = c_1( T_{\fol}^*)c_1([D]^*)$ and we get
$$
\displaystyle \sum^{n}_{j=0}\left(\sum^{n -j}_{i=0}c_{n-j-i}(TX)c_1([D]^*)^i\right)\hspace{-0.1cm} c_1( T_{\fol}^*)^j =\hspace{-0.1cm} \displaystyle\sum^{n}_{i=0}\hspace{-0.1cm}\left(\sum^{n -i}_{j=0}c_{n-i-j}(TX)c_1( T_{\fol}^*)^j\right)\hspace{-0.1cm} c_1([D]^*)^i.
$$

\noindent Thus we have that

{\small{
\begin{eqnarray}\nonumber
\hspace{-8.2cm } \displaystyle\int_{X}\hspace*{-0.15cm }c_{n}\hspace*{-0.35cm }&(&\hspace*{-0.35cm }TX - [D]- T_{\fol}) 
\\\nonumber && \\ \nonumber =&& \hspace*{-0.9cm }\displaystyle\int_{X}\displaystyle\sum^{n}_{i=0}\left(\sum^{n -i}_{j=0}c_{n-i-j}(TX)c_1( T_{\fol}^*)^j\right)c_1([D]^*)^i
\\\nonumber && \\ \nonumber = &&\hspace*{-0.9cm } \displaystyle\int_{X} \hspace*{-0.1cm }\displaystyle\sum^{n}_{j=0}c_{n-j}(TX)c_1(T_{\fol}^*)^j + \displaystyle\int_{X}\displaystyle\sum^{n-1}_{i=1}\left(\sum^{n -i}_{j=0}c_{n-i-j}(TX)c_1( T_{\fol}^*)^j\right)c_1([D]^*)^i  + \int_{X}c_1([D]^{\ast})^n  \\\nonumber && \\ \nonumber = && \hspace*{-0.9cm } \displaystyle\int_{X}\hspace*{-0.2cm }\displaystyle c_{n}(TX - T_{\fol}) \hspace*{-0.07cm } + \hspace*{-0.14cm } \displaystyle\int_{X}\hspace*{-0.07cm }\displaystyle\sum^{n-1}_{i=1}\hspace*{-0.09cm }\left(\hspace*{-0.07cm }c_{n-i}(TX) \hspace*{-0.07cm } + \hspace*{-0.14cm }  \sum^{n -i}_{j=1}c_{n-i-j}(TX)c_1( T_{\fol}^*)^j\hspace*{-0.09cm }\right)\hspace*{-0.09cm }c_1([D]^*)^i  \hspace*{-0.07cm } + \hspace*{-0.14cm } \int_{X}c_1([D]^{\ast})^n  \\\nonumber && \\ \nonumber = && \hspace*{-0.9cm } \displaystyle\int_{X}\hspace*{-0.2cm }\displaystyle c_{n}(TX - T_{\fol}) \hspace*{-0.07cm } + \hspace*{-0.14cm }\displaystyle\int_{X}\displaystyle\sum^{n}_{i=1}c_{n-i}(TX)c_1([D]^*)^i \hspace*{-0.07cm } + \hspace*{-0.14cm } \displaystyle\int_{X}\hspace*{-0.05cm }\displaystyle\sum^{n-1}_{i=1}\hspace*{-0.11cm } \left( \sum^{n -i}_{j=1}c_{n-i-j}(TX)c_1( T_{\fol}^*)^j\hspace*{-0.11cm }\right)\hspace*{-0.10cm }c_1([D]^*)^i  
\\\nonumber && \\ \nonumber = && \hspace*{-0.9cm } \displaystyle\int_{X}\hspace*{-0.18cm }\displaystyle c_{n}(TX - T_{\fol}) \hspace*{-0.07cm } + \hspace*{-0.14cm }\displaystyle\int_{X}\displaystyle \hspace*{-0.26cm } c_{n-1}(TX - [D])c_1([D]^*) \hspace*{-0.07cm } + \hspace*{-0.14cm } \displaystyle\int_{X}\hspace*{-0.10cm }\displaystyle\sum^{n-1}_{j=1}\hspace*{-0.10cm } \left( \sum^{n -j}_{i=1}\hspace*{-0.05cm }c_{n-j-i}(TX)c_1([D]^*)^i\right)c_1( T_{\fol}^*)^j \\\nonumber && \\ \nonumber = && \hspace*{-0.90 cm } \displaystyle\int_{X}\hspace*{-0.29cm }\displaystyle c_{n}\hspace*{-0.04cm }(\hspace*{-0.04cm }TX - T_{\fol}\hspace*{-0.04cm }) \hspace*{-0.07cm } - \hspace*{-0.16cm }
\displaystyle \int_{X}\hspace*{-0.27cm }\displaystyle c_{n-1}\hspace*{-0.04cm }(\hspace*{-0.035cm }TX - [D])c_1([D]) \hspace*{-0.07cm } - \hspace*{-0.14cm } \displaystyle\int_{X}\hspace*{-0.1cm }\displaystyle\sum^{n-1}_{j=1}\hspace*{-0.1cm }\left( \sum^{n -j}_{i=1}\hspace*{-0.03cm }c_{n-j-i}(TX)(-1)^{i-1}\hspace*{-0.1cm }c_1([D])^i\hspace*{-0.1cm }\right)\hspace*{-0.1cm }c_1( T_{\fol}^*)^j   \\\nonumber &&  \\ \nonumber =&& \hspace*{-0.9cm } \sum_{p\in Sing(\fol)}\hspace{-0.40cm}\mu_p(\fol) \hspace*{-0.07cm } - \hspace*{-0.14cm } \displaystyle \int_{D}\displaystyle c_{n-1}(TX - [D]) \hspace*{-0.07cm } - \hspace*{-0.14cm } \displaystyle\int_{D}\displaystyle\sum^{n-1}_{j=1}\left( \sum^{n -j}_{i=1}c_{n-j-i}(TX)(-1)^{i-1}c_1([D])^{i-1}\right)c_1( T_{\fol}^*)^j  
\end{eqnarray}
}}

\noindent where in the last equality we have used the Baum-Bott classical formula (\ref{bbformula}) and the fact that $c_1([D])$ is Poincaré dual to
the fundamental class of $D$. On the other hand, using (\ref{gaf}), we obtain
{\small{
\begin{eqnarray}\nonumber
 & & \hspace*{-0.95 cm} \displaystyle\int_{X}c_{n}(TX - [D]- T_{\fol})
=\\\nonumber && \\ \nonumber =& & \hspace*{-0.95 cm}\sum_{p\in Sing(\fol)}\hspace{-0.45cm}\mu_p(\fol) \hspace*{-0.07cm } - \hspace*{-0.07cm } \displaystyle \chi(D) \hspace*{-0.07cm } + \hspace*{-0.07cm }(-1)^n\hspace{-0.35cm}\sum_{p\in Sing(D)}\hspace{-0.45cm} \mu_p(D)  \hspace*{-0.07cm } - \hspace*{-0.14cm } \displaystyle\int_{D}\hspace{-0.09cm}\displaystyle\sum^{n-1}_{j=1}\hspace{-0.1cm}\left( \sum^{n -j}_{i=1}\hspace{-0.1cm}c_{n-j-i}(TX)(-1)^{i-1}\hspace{-0.1cm}c_1([D])^{i-1}\hspace{-0.1cm}\right)\hspace{-0.1cm}c_1( T_{\fol}^*)^j, 
\end{eqnarray}
}}
\noindent and, we compute the last integral as follows:
{{
\begin{eqnarray}\nonumber
&& \hspace{-0.65 cm} \displaystyle\int_{D}\displaystyle\sum^{n-1}_{j=1}\left( \sum^{n -j}_{i=1}c_{n-j-i}(TX)(-1)^{i-1}c_1([D])^{i-1}\right)c_1( T_{\fol}^*)^j \\\nonumber && \\ \nonumber &= & \displaystyle\int_{D}\displaystyle\sum^{n-1}_{j=1}c_{n-1-j}((TX)-[D])c_1( T_{\fol}^*)^j \\\nonumber && \\ \nonumber &=& \displaystyle\int_{D}\hspace{-0.25 cm} c_{n-1}(TX - [D]- T_{\fol}) -\int_{D}\hspace{-0.25 cm} c_{n-1}(TX - [D]) \\\nonumber && \\ \nonumber &=& \displaystyle - \chi(D) + \hspace{-0.25 cm} \sum_{p\in  S(\fol, D) }\hspace{-0.4 cm} Sch_p(\fol,D)
\end{eqnarray}
}}

\noindent where in the last equality we have used the formula (\ref{fof}). This proves the formula of item (III).

The item (IV) follows directly of item (III). In any case, we sketch the argument for the convenience of the reader: considering the following decompositions into disjoint unions
\begin{eqnarray}\nonumber
Sing(\fol) &=& (Sing(\fol) \cap (X - D)) \cup (Sing(\fol)\cap D_{reg}) \cup (Sing(\fol) \cap Sing(D))\\\nonumber
Sing(D) &=& (Sing(\fol) \cap Sing(D)) \cup (Sing(D) - Sing(\fol))
\\\nonumber
S(\fol,D) &=& (Sing(\fol) \cap D_{reg}) \cup (Sing(\fol) \cap Sing(D)) \cup (Sing(D) - Sing(\fol))  
\end{eqnarray}

\noindent we can write the formula in (III) in the following format
{ \small{
\begin{eqnarray}\label{fgh}
&&\hspace*{-1cm}\int_{X} \hspace{-0.10 cm} c_n (TX - [D] - T_{\fol}) \hspace*{0.05cm} = \hspace*{-0.25cm}\sum_{p\in Sing(\fol)\cap (X- D)}\hspace{-0.7 cm} \mu_{p}(\fol) \,\,\,+ \hspace*{-0.25cm}\sum_{p\in Sing(\fol)\cap D_{reg}}\hspace{-0.7 cm} (\mu_{p}(\fol) - Sch_p(\fol,D))\\\nonumber  && \\\nonumber  && \hspace*{-0.95cm} + \hspace*{-0.25cm}\sum_{p\in Sing(\fol)\cap D_{reg}}\hspace{-0.35 cm} (\mu_{p}(\fol) - Sch_p(\fol,D) + (-1)^n\mu_p(D)) \\\nonumber  && \\\nonumber  && \hspace*{-0.95cm} + \hspace*{-0.25cm}\sum_{p\in Sing(D)- Sing(\fol)}\hspace{-0.35 cm} (- Sch_p(\fol,D) + (-1)^n\mu_p(D)).
\end{eqnarray}
}}

\noindent and the item (IV)
follows from ($\ref{fgh}$).

By using item (III), we obtain

\begin{eqnarray}\nonumber
\displaystyle\sum_{p\in  S(\fol, D) }\hspace*{-0.3cm} Sch_p(\fol,D) \hspace*{-0.25cm} &=& \hspace*{-0.45cm}\sum_{p\in  Sing(\fol) } \hspace{-0.35 cm} \mu_p(\fol) - \int_{\PP^n}\hspace*{-0.2cm} c_n (T\PP^n - [D] - T_{\fol})\,\, + \,\, (-1)^{n}\hspace{-0.25 cm}\sum_{p\in  Sing(D) } \hspace{-0.35 cm} \mu_p(D)\\\nonumber
&=& \sum_{i=0}^{n-1}\left(1 - (1-k)^{n-i}\right)d^i \,\,\,\, + \,\,\,\,(-1)^n\hspace{-0.45 cm}\displaystyle\sum_{p\in Sing(D)}\mu_p(D),
\end{eqnarray}

\noindent where in the last step we have used the formula (\ref{2121}) and Lemma 3.1.

\noindent $\square$

\begin{lema} \label{lemx}
Let $\fol$ be a foliation of dimension one on $\PP^n$ and degree $\deg(\fol) = d$ and let $D \subset \PP^n$ be an irreducible hypersurface of degree $\deg(D) = k$ with isolated singularities and such that $\fol$ leaves $D$ invariant. Then,
\begin{eqnarray}\nonumber\label{hhhh}
\displaystyle \int_{\PP^n}c_{n}(T\PP^n - [D]- T_{\fol}) =  \displaystyle\sum^{n}_{i=0}(-1)^i(k-1)^id^{n-i}.
\end{eqnarray}
\end{lema}

\noindent {\bf Proof.:} According to \cite[Lemma 4.1]{DSM}, the top Chern number of the virtual bundle $T\PP^n - [D] - T\fol$ is given by

\begin{eqnarray}\nonumber\label{ff1ff}
\displaystyle \int_{\PP^n}c_{n}(T\PP^n - [D]- T_{\fol}) =  \displaystyle\sum^{n}_{i=0} \sum^{n-i}_{j=0}\binom{n + 1}{n-i-j}
(-k)^j (d-1)^i.
\end{eqnarray}

\noindent On the other hand, by \cite[Lemma 4.2]{MD}, the expression 

$$\displaystyle\sum^{n}_{i=0} \sum^{n-i}_{j=0}\binom{n + 1}{n-i-j}
(-k)^j (d-1)^i$$ 

\noindent is equivalent to
$$
\sum_{i=0}^n(-1)^i(k-1)^id^{n-i}.
$$ 

\noindent This proves the Lemma. 

\noindent $\square$
\bigskip

\hyphenation{po-si-ti-ve}

\section{Proof of Theorem 2}\label{a4}
Suppose $p \in Sing(\fol) \cap Sing(D)$ (if $p\in D_{reg}$,  then the GSV index and the Schwartz index are both equal to the Poincaré-Hopf index, and therefore are positive) and  let $\{f=0\}$ be a local equation of $D$ in a neighborhood $U$ of $p$. By taking a smaller U, if necessary, let $v = \sum_{i=1}^na_i\frac{\partial}{\partial z_i}$ be the holomorphic vector field inducing the foliation $\fol$ on $U$, with holomorphic functions $a_i\in \OO(U)$.  Since the GSV index coincides with the {\it homological index} (see for example \cite{ABC},\cite{a08}), it follows from \cite[Theorem 1]{a08} that:

\medskip

\noindent For even $n$:
{\small{
\begin{eqnarray}\label{eq8888}
 GSV_p(\fol, D) = \displaystyle dim_{\CC}\frac{\OO_{n,p}}{\left\langle f,a_1,\ldots, a_n\right\rangle} - dim_{\CC}\frac{\OO_{n,p}}{\left\langle f,J_f\right\rangle}.
\end{eqnarray}
}} 

\noindent For odd $n$:

{\small{
\begin{eqnarray}\label{eq888}
 \hspace*{1cm} GSV_p(\fol, D) =  \displaystyle dim_{\CC}\frac{\OO_{n,p}}{\left\langle a_1,\ldots, a_n\right\rangle} -  dim_{\CC}\frac{\OO_{n,p}}{\left\langle \displaystyle \frac{df}{f}(v), a_1,\ldots, a_n\right\rangle} + dim_{\CC}\frac{\OO_{n,p}}{\displaystyle\left\langle f,J_f\right\rangle}.
\end{eqnarray}
}}

\noindent Here $\OO_{n,p}$ denotes the ring of germs of holomorphic functions at $(X, p)$ and $J_f = \left\langle \frac{\partial f}{\partial z_1},\ldots,\frac{\partial f}{\partial z_n}\right\rangle \subset \OO_{n,p}$ denotes the Jacobian ideal of $f$. 

We observe that 
\begin{eqnarray}\nonumber\label{ne}
\displaystyle dim_{\CC}\frac{\OO_{n,p}}{\left\langle f,a_1,\ldots, a_n\right\rangle} > 0
\end{eqnarray}
\noindent and
\begin{eqnarray}\nonumber\label{me}
\displaystyle dim_{\CC}\frac{\OO_{n,p}}{\left\langle f,J_f\right\rangle} > 0
\end{eqnarray}

\noindent (otherwise, there exist germs $g_i's$ in $\OO_{n,p}$ such that  $1 = \sum_{i=1}^ng_ia_i + g_{n+1}f$ or $1 = \sum_{i=1}^ng_i\frac{\partial f}{\partial z_i} + g_{n+1}f$ , which is a contradiction to the fact $p\in Sing(\fol) \cap Sing(D)$).

Since $\left\langle f,J_f\right\rangle \OO_{n,p}\supset \left\langle J_f\right\rangle\OO_{n,p}$, it follows from (\ref{form1})  and (\ref{eq8888}) that
\begin{eqnarray}\nonumber
0 < \displaystyle dim_{\CC}\frac{\OO_{n,p}}{\left\langle f,a_1,\ldots, a_n\right\rangle} \,\, \leq\,\,\ GSV_p(\fol, D) \,\, + \,\, \mu_p(D) \,\, = \,\, Sch_p(\fol, D).
\end{eqnarray}

On the other hand, since $\left\langle  \frac{df}{f}(v), a_1,\ldots, a_n\right\rangle\OO_{n,p} \supset \left\langle \displaystyle a_1,\ldots, a_n\right\rangle\OO_{n,p}$ it follows from (\ref{eq888}) that
\begin{eqnarray}\nonumber
0 < dim_{\CC}\frac{\OO_{n,p}}{\left\langle f,J_f\right\rangle} \,\, \leq\,\,\ GSV_p(\fol, D).
\end{eqnarray}
$\square$

The number $dim_{\CC}\frac{\OO_{n,p}}{\left\langle f,J_f\right\rangle}$ that appears in the above demonstration is the {\it Tjurina number} $\tau_p(D)$ of $D$ at $p$ (see for example, \cite{sai000}). Therefore, from Theorem 2 we deduce that the Tjurina number is a lower bound for GSV index, in odd  dimension. Furthermore, in the case where $X$ is the $n$-dimensional projective space $\PP^n$, using the \cite[Corolary 2.1]{lll}, we obtain a lower bound in terms of the (local) multiplicity of the hypersurface $D$:

\begin{cor}\label{coro} Let $X$ be an $n$-dimensional complex manifold, $D\subset X$ a hypersurface with isolated singularities and $\fol$ a one-dimensional holomorphic foliation on $X$, with isolated singularities, which leaves $D$ invariant. Suppose $n$ odd. If $p\in Sing(\fol)\cap Sing(D)$ then one has $\tau_p(D) \leq GSV_p(\fol,D)$. In particular, if $X = \PP^n$ and $m$ denotes the multiplicity of $D$ at $p$, one has $\displaystyle\frac{(m-1)^n}{n} \leq GSV_p(\fol,D)$.
\end{cor}

\medskip

\section{Proof of Theorem 3}
To prove Theorem 3 we consider first the following lemma

\begin{lema}\label{tta}
Let $k, d$ and $n$ natural numbers, with $k\geq 1$, $d \geq 0$ and $n\geq 2$. If $k>d+2$, then 
\begin{eqnarray}\nonumber
\sum_{i=0}^{n-1}\left(1 - (1-k)^{n-i}\right)d^i < 0.
\end{eqnarray}
\end{lema}
\noindent {\bf Proof.:} The original sum can be written as follows
\begin{eqnarray}\nonumber
\sum_{i=0}^{n-1}\left(1 - (1-k)^{n-i}\right)d^i &=&  \sum_{i=0}^{n-1}d^i \hspace{ 0.15cm} + \displaystyle\sum_{\mbox{{\footnotesize $j$=$1$,$3$,...,$n$$-$$1$}}}\hspace{- 0.55cm}(k-1)^{n-i}d^i \hspace{ 0.15cm} - \hspace{- 0.12cm}\displaystyle\sum_{\mbox{{\footnotesize $j$=$0$,$2$,...,$n$$-$$2$}}}\hspace{- 0.5cm}(k-1)^{n-i}d^i \hspace{- 0.08cm}.
\end{eqnarray}

\noindent By using the hypothesis $k-1>d+1$ in the last sum, we obtain
\begin{eqnarray}\noindent \nonumber
\hspace{- 0.15cm} -\hspace{- 0.15cm} \displaystyle\sum_{\mbox{{\footnotesize $j$=$0$,$2$,...,$n$$-$$2$}}}\hspace{- 0.5cm}(k-1)^{n-i}d^i \hspace{- 0.08cm} &<& - \displaystyle\sum_{\mbox{{\footnotesize $j$=$0$,$2$,...,$n$$-$$2$}}}\hspace{- 0.5cm}(d+1)(k-1)^{n-i-1}d^i \hspace{- 0.08cm} \\\nonumber && \\\nonumber
&=& - \displaystyle\sum_{\mbox{{\footnotesize $j$=$0$,$2$,...,$n$$-$$2$}}}\hspace{- 0.5cm}(k-1)^{n-i-1}d^{i+1} \hspace{ 0.4cm} -  \displaystyle\sum_{\mbox{{\footnotesize $j$=$0$,$2$,...,$n$$-$$2$}}}\hspace{- 0.58cm}(k-1)^{n-i-1}d^{i} \\\nonumber && \\\nonumber
&<& - \displaystyle\sum_{\mbox{{\footnotesize $j$=$0$,$2$,...,$n$$-$$2$}}}\hspace{- 0.5cm}(k-1)^{n-i-1}d^{i+1} \hspace{ 0.4cm} -  \displaystyle\sum_{\mbox{{\footnotesize $j$=$0$,$2$,...,$n$$-$$2$}}}\hspace{- 0.5cm}(d+1)d^{i} \\\nonumber && \\\nonumber
&=& - \displaystyle\sum_{\mbox{{\footnotesize $j$=$1$,$3$,...,$n$$-$$1$}}}\hspace{- 0.5cm}(k-1)^{n-i}d^i \hspace{ 1.1cm}  -  \hspace{ 0.4cm}\sum_{i=0}^{n-1}d^i,\end{eqnarray}
\noindent and the desired inequality is proved.

\noindent $\square$

\noindent {\it \bf Prove of Theorem 3.}
Let $deg(D) = k$, and let $deg(\fol) = d$. By item (V) of Theorem 1, we obtain 
\begin{eqnarray}\nonumber
\hspace{0.5cm}\displaystyle\sum_{p\in  S(\fol, D) } \hspace{-0.4cm} Sch_p(\fol, D) &=&  \sum_{i=0}^{n-1}\left(1 - (1-k)^{n-i}\right)d^i \,\,\,\, + \,\,\,\,\hspace{-0.45 cm}\displaystyle\sum_{p\in Sing(D)}\mu_p(D) \Longleftrightarrow \\\nonumber && \\\nonumber \sum_{i=0}^{n-1}\left(1 - (1-k)^{n-i}\right)d^i &=& 
\displaystyle\sum_{p\in  Sing(D) } \hspace{-0.4cm}(Sch_p(\fol, D) - \mu_p(D)) + \sum_{p\in  Sing(\fol)\cap D_{reg} } \hspace{-0.8cm} Sch_p(\fol, D).    
\end{eqnarray}
\noindent  Since Schwartz index is positive and $\mu_p(D) =1$, for all $p\in Sing(D)$, we get $\displaystyle \sum_{p\in  Sing(D)} \hspace{-0.4cm} (Sch_p(\fol, D) - \mu_p(D))\geq 0$ and $\displaystyle\sum_{p\in  Sing(\fol)\cap D_{reg} } \hspace{-0.85cm} Sch_p(\fol, D)\geq 0$, and consequently

$$
\sum_{i=0}^{n-1}\left(1 - (1-k)^{n-i}\right)d^i \geq 0.
$$
\noindent Therefore, by Lemma \ref{tta}, we obtain

$$
k \leq d+2.
$$
$\square$

\section{Proof of Theorem 4}

Let $deg(D) = k$, and let $deg(\fol) = d$. It follows from the item (V) of Theorem 1 that
$$
\sum_{i=0}^{n}(1-k)^{n-i}d^i =  \sum_{i=0}^{n}d^i + \displaystyle\sum_{p\in  Sing(D) } \hspace{-0.4cm}(\mu_p(D) - Sch_p(\fol, D)) \hspace{0.2cm} - \hspace{-0.2cm} \sum_{p\in  Sing(\fol)\cap D_{reg} } \hspace{-0.8cm} Sch_p(\fol, D).  
$$

\noindent Since Schwartz index is positive, we obtain
$$
\sum_{i=0}^{n}(1-k)^{n-i}d^i \,\,\,\,\leq\,\,\,\,  \sum_{i=0}^{n}d^i + \displaystyle\sum_{p\in  Sing(D) } \hspace{-0.4cm}(\mu_p(D) -1).
$$
Thus, using the relation
\begin{eqnarray}\nonumber
\sum_{j=0}^nd^j  - \sum_{i=0}^{n}(1-k)^{n-i}d^i  =  \displaystyle \displaystyle  \sum^{n-1}_{i=1}\sum^{i-1}_{j=0}\binom{i}{j}(-1)^{i-j+1}k^{i-j}d^{n-i} -  \sum^{n-1}_{j=0}\binom{n}{j}(-k)^{n-j}      
\end{eqnarray}
\noindent we get the desired inequality
\begin{eqnarray}\nonumber
\displaystyle \sum^{n-1}_{j=0}\binom{n}{j}(-k)^{n-j} - \hspace*{-0.45cm}\sum_{p\in Sing(D)}\hspace*{-0.25cm} (\mu_p(D)- 1) \leq \displaystyle \sum^{n-1}_{i=1}\sum^{i-1}_{j=0}\binom{i}{j}(-1)^{i-j+1}k^{i-j}d^{n-i}.
\end{eqnarray} 
$\square$

\begin{exe}\label{exeexe}
Let $\PP^4$ be the complex projective space of dimension 4, with homogeneous coordinates $(X_0:X_1:X_2:X_3:X_{4})$, and let $D\subset \PP^4$ be the hypersurface defined by 
$$
X_1^k + X_2^k +X_3^k+ X_4^k = 0,
$$ 
\noindent with $k>1$, whose singular set is $Sing(D) = \{(1:0:0:0:0)\}$. 

In the affine space
$X_0 \neq 0$, with coordinates 
$$
\displaystyle (x_1,x_2,x_3,x_4) = \left(X_1/X_0,X_2/X_0,X_3/X_0,X_{4}/X_0 \right),
$$ 
\noindent the defining function of $D$ is given  by 
$$
f(x_1,x_2,x_3,x_4) = x_1^k + x_2^k +x_3^k +x_4^k.
$$  
\noindent Since $f$ a {\it Pham-Brieskorn} polynomial, the Milnor number of $D$ in singular point $p_0=(1:0:0:0:0)$  is given by (see \cite{JM})
$$
\mu_{p_0}(D) = (k-1)^4.
$$
\noindent On the other hand, the foliation $\fol$ on $\PP^4$ defined by the vector field (affine coordinate $X_0 \neq 0$)
$$
v =  x_1\displaystyle\frac{\partial}{\partial x_1} + x_2\displaystyle\frac{\partial}{\partial x_2} + x_3\displaystyle\frac{\partial}{\partial x_3} + x_4\displaystyle\frac{\partial}{\partial x_4}
$$
\noindent  leaves D invariant and the singular set of $\fol$ consist of point $p_0 = (1:0:0:0:0)$. 
\medskip

\noindent {\bf(i)} We have 
\begin{eqnarray}\nonumber
\sum_{i=0}^3(1-(1-k)^{4-i})d^i + \sum_{p\in Sing(D)}\mu_p(D) &=& 1-(1-k)^{4} + (1-k)^{4}\\\nonumber  &=& 1 
\end{eqnarray}
\noindent and, since $v$ a radial vector field,
$$
\sum_{p\in S(\fol,D)} Sch_p(\fol,D) = Sch_{p_0}(\fol,D) = 1,
$$ 
\noindent according to item (V) of Theorem 1.

\medskip

\noindent {\bf(ii)} If $k=2$, then $p_0$ have multiplicity $1$ and  $deg(D) \leq  deg(\fol) + 2$  according to Theorem 3.
\end{exe}

\section{Proof of Theorem 5}\label{b0}
 
Let $D$ be a complex compact hypersurface with isolated singularities $p_1,\ldots,p_{s_1}$ in a complex manifold $X$ of even dimension. Suppose that $D$ admits a holomorphic vector field singular at the $p_j's$. By using \cite[Theorem 2.1.1]{BarSaeSuw}, we obtain
$$
\chi (D) = \sum_{j=1}^{s_1} Sch_{p_j}(\fol, D) + \sum_{j=1}^{s_2} Sch_{q_j}(\fol, D)> \sum_{j=1}^{s_1} Sch_{p_j}(\fol, D) + s_2,
$$
\noindent where in the last inequality we have used the fact that at regular points of $D$ the Schwartz index is equal to the Poincaré-Hopf index, and therefore is positive. Furthermore, since $dim(X)$ is even, then it follows from item (i) of Theorem 2 that $Sch_{p_j}(\fol, D)>0$ for each $j = 1,\ldots,s_1$ and thus $\sum_{j=1}^{s_1} Sch_{p_j}(\fol, D)>s_1$. 

\noindent $\square$

\end{document}